\documentclass[final,3p,times]{elsarticle}
\usepackage{graphics}
\usepackage[utf8]{inputenc} 
\usepackage{amssymb,amsthm,amsmath,amsfonts}
\journal{}
\usepackage{url}
  \newtheorem{theorem}{Theorem}
  \newtheorem{proposition}{Proposition}
               \newtheorem{corollary}{Corollary}
							
               \newtheorem{lemma}{Lemma}
               \newtheorem{example}{Example}

               \def\pf{\par\noindent {\em Proof.}~\par\noindent}
               \def\qed{~\hfill{$\square$}\pagebreak[1]\par\medskip\par}

\newcommand{\f}{\bf{infra}}
\newcommand{\cD}{{\cal D}}

\newcommand{\cP}{{\cal P}}

\newcommand{\mdiv}{{\mbox{div}}}

\newcommand{\R}{{\mathbb R}}
\newcommand{\B}{{\mathbb B}}

\newcommand{\I}{{\cal I}}

\newcommand{\grad}{{\rm grad}}

\newcommand{\curl}{{\rm curl}}
\newcommand{\ux}{\underline{x}}

\newcommand{\uy}{\underline{y}}

\newcommand{\pux}{\partial_{\ux}}

\newcommand{\pD}{{^\vartheta\!\pux}}
\newcommand{\hD}{{^\chi\!\pux}}

\newcommand{\cS}{{\cal S}}

\newcommand{\E}{{\bf E}}

\begin{document}

\begin{frontmatter}

\title{The Dirichlet problem for a non-commutative elliptic operator in the ball}

\author{Arsenio Moreno Garc\'ia$^{(1)}$; Daniel Alfonso Santiesteban$^{(2)*}$; Ricardo Abreu Blaya$^{(2)}$}
\address{
$^1$Faculty of Informatics and Mathematics, University of Holguín, Cuba\\
$^2$Faculty of Mathematics, Autonomous University of Guerrero, Mexico}
\ead{amorenog@uho.edu.cu, danielalfonso950105@gmail.com, rabreublaya@yahoo.es}

\cortext[cor1]{Corresponding author}

\begin{abstract}
In this paper we establish a necessary and sufficient condition for the solvability of the Dirichlet problem in the unit ball $\B$ for the second order system $\pux f\pux=0$, where $\pux$ stands for the Clifford-algebra valued Dirac operator in $\R^m$. We prove that the problem admits a solution whenever the boundary data belong to $C^1(\partial\B)$. Conversely, we construct explicit counterexamples showing that solvability may fail if this smoothness assumption on the boundary is not satisfied. Moreover, in contrast with the standard commutative setting, we show how that solutions exhibit a regularity one order less than the boundary data, when the latter belong to $C^{k,\nu}(\partial\B)$.
\end{abstract}

\begin{keyword}
Dirichlet problem, Clifford algebras, Dirac operator, inframonogenic functions.\\
\noindent\textit{2020 Mathematics Subject Classification: 30G35, 35J47, 35J57.}
\end{keyword}

\end{frontmatter}
\section{Introduction}
We consider in $\R^m$ the Dirac operator
\[
\pux:=e_1{\partial_{x_1}}+e_2{\partial_{x_2}}+\cdots +e_m{\partial_{x_m}},
\]
where $e_1,e_2,\dots ,e_m$ stand for orthonormal basis vectors of $\R^m$, which subjected to the multiplication rules
\[
e_i^2=-1,\quad e_ie_{j}=-e_{j}e_i,\quad i,j=1,2,\dots,m,\quad i<j,
\]
generate a minimal enlargement of $\R^m$  to a real linear associative (but not commutative) algebra, namely the Clifford algebra $\R_{0,m}$.

The Dirac operator factorizes the Laplacian in the sense that $\pux\pux=-\Delta$. Thus it is seen that the Laplace equation may be written as 
\begin{equation}\label{L}
\pux\pux f=0.
\end{equation}
Inframonogenic functions arise from a non-commutative version of \eqref{L}, i.e., as the solutions of the sandwich equation 
\begin{equation}\label{I}
\pux f\pux=0.
\end{equation}
Such functions were originally introduced in \cite{MPS1,MPS2} and more recently in \cite{MAB1,MAB2,MAB3} have been found interesting connections between them and the solutions of the Lam\'e-Navier system in linear elasticity theory. 

It should be noted here that when restricting our attention to vector-valued solutions $f=\vec{u}$ in $\R^3$, the above sandwich equation \eqref{I} can be written in the form
\[
\grad\,\mdiv\,\vec{u}+\curl^2\,\vec{u}=0,
\] 
which has already been discussed in \cite[pag. 267]{Dzhu} as an example of elliptic but not strongly elliptic equation. This system arises in several areas of mathematics and physics. In elasticity theory, it appears in the Navier–Lamé equations governing displacement fields in solids (see for example \cite{Achenbach1973, Mura1987,Ciarlet1988}). In fluid dynamics, related formulations are found in the study of incompressible viscous flows \cite{Batchelor1967,Ladyzhenskaya1969}. 
Mathematically, it is directly linked to the vector calculus identity underlying the Helmholtz decomposition of vector fields \cite{Arfken2005, Cantarella2000}.

In this paper we shall focus our attention on the solvability of the Dirichlet problem 
\begin{equation}\label{DPI}
\Biggl\{
\begin{array}{rl}
\pux f\pux=0 &\,\,\mbox{in}\,\,\B,\\f=g&\,\,\mbox{in}\,\,\partial\B,
\end{array} 
\end{equation}
where $\B$ is the unit ball in $\R^m$ and $g$ is assumed to be continuous on $\partial\B$.
   
\section{Preliminaries}
This section summarizes a series of preliminary notions and tools that are required in the main part of this work.

\subsection{Algebraic preliminaries}
 
As mentioned in the introduction the Clifford algebra $\R_{0,m}$, in which the Euclidean space $\R^m$ is embedded, is generated by the orthonormal basis $e_1,e_2,\dots, e_m$. Thus we describe the elements of $\R_{0,m}$ in the form $a=\sum_{A} a_A e_A$, where as indices we use the elements $A$ of the set containing the ordered subsets of $\{1,2,\ldots ,m\}$, with the empty subset corresponding to the index $0$.

In particular, ${Sc}[a]:=a_0$ represents the scalar part of $a$. The mapping $a\mapsto\overline{a}$, with $\overline{e_i}=-e_i$ is called conjugation in $\R_{0,m}$. The norm of an element $a\in\R_{0,m}$ is defined by $\|a\|^2=Sc[a\overline{a}]$. When $\ux\in\R^m$, the norm of $\ux$ is obviously equal to the Euclidean norm $|\ux|$.

The linear subspace of $\R_{0,m}$ spanned by the $\binom{m}{p}$ products $e_A$, $|A|=p$, is denoted by $\R_{0,m}^{(p)}$. The elements of $\R_{0,m}^{(p)}$ are the so-called $p$-vectors.

Let $[\,]_{p}:\R_{0,m}\mapsto\R_{0,m}^{(p)}$ be a linear projection with 
\[
[a]_p=\sum_{|A|=p}a_A\,e_A.
\] 
In this way an arbitrary element $a\in\R_{0,m}$ can be written in the form
\[
a=\sum_{p=0}^m[a]_{p}.
\]
We will make repeated use of the operator $\Psi:\R_{0,m}\mapsto\R_{0,m}$ given by
\[
\Psi(a)=\sum_{j=1}^m e_j a e_j,
\]
for $a\in\R_{0,m}$.

We easily check that 
\begin{equation}\label{=0}
\Psi(Y_p)=(-1)^{p+1}(m-2p)Y_p,
\end{equation}
when $Y_p$ is a $p$-vector (see \cite{MPS2}).

Moreover, the operator $\Psi$ is a bijective map when dealing with odd dimension, and its inverse mapping is given by (see \cite[Proposition 1]{M2A})
\[
\Psi^{-1}(a)=\sum_{p=0}^m\frac{(-1)^{p+1}}{m-2p}[a]_p.
\] 
\subsection{Analytic preliminaries}
We will consider $\R_{0,m}$-valued functions on subsets of $\R^{m}$, which can be expressed as $f=\sum_{A} f_A e_A$, with $f_A$ being $\R$-valued. In this context the spaces of $k$-times H\"older continuous differentiable functions have the usual component-wise meaning and will be denoted by $C^{k,\alpha}(\E)$, where $\E$ is a subset of $\R^{m}$.

In particular, an $\R_{0,m}$-valued function $f$ is assumed to be in the $C^{1,\nu}(\E)$, $0<\nu<1$, if each real component $f_A$  of $f$ does so.
  
Let us come back to the previously defined Dirac operator $\pux$.

An $\R_{0,m}$-valued function $f$ is called left monogenic (right monogenic) in a domain $\Omega\subset\R^m$ if $\pux f = 0$ ($f\pux= 0$) in $\Omega$ (see for example \cite{BDS}).  More generally, left (right) polymonogenic functions of order $k$ are the $\R_{0,m}$-valued solutions of the iterated Dirac equation $\pux^k\,f = 0$ ($f\pux^k = 0$). These functions are considered for example in \cite{Br1,Br2,Ry}. In this paper we make particular use of polymonogenic functions of order $3$, also referred to as three-monogenic functions. Two-sided monogenic (polymonogenic) functions are those that are simultaneously left and right monogenic (polymonogenic).

Let $F_k$ be a $k$-vector valued function (or $k$-vector field). The action of the Dirac operator\index{Dirac!operator} on $F_k$ can be decomposed as
\begin{equation}
\pux F_k=\pux\bullet F_k+\pux\wedge F_k,
\end{equation}
where
\begin{equation}
\pux\bullet F_k=\frac{1}{2}[\pux F_k-(-1)^kF_k\pux]=[\pux F_k]_{k-1}
\end{equation}
and
\begin{equation}
\pux\wedge F_k=\frac{1}{2}[\pux F_k+(-1)^kF_k\pux]=[\pux F_k]_{k+1}.
\end{equation}
The inner ($\bullet$) and outer ($\wedge$) products\index{inner product}\index{outer product} play a central role in the reformulation and factorization of the several equations in Clifford analysis. The sandwich operator and the Laplacian\index{Laplace!operator} preserve the space of $k$-vector valued functions. This follows from the identities
\begin{align}
-\Delta F_k&=\pux\cdot\pux\wedge F_k+\pux\wedge\pux\cdot F_k,\label{deltarel}\\
\pux F_k\pux&=(-1)^k(\pux\cdot\pux\wedge F_k-\pux\wedge\pux\cdot F_k),\label{sandrel}
\end{align}
for $F_k\in C^2(\Omega,\R_{0,m}^{(k)})$. In particular, \eqref{deltarel} and \eqref{sandrel} show that the sandwich operator $\pux(\cdot)\pux$ behaves analogously to the Laplacian when acting on scalar and pseudoscalar fields\index{pseudoscalar}.

Finally, to connect the Dirac operator $\pux$ to the previously defined operator $\Psi$, the following formulas will be useful (see \cite[Lemma 2.1]{MAB3}). 
\begin{proposition}\label{p1}
Let be $f$ be twice continuously differentiable, then
\begin{itemize}
\item[(i)] $\pux (f\,\ux)=(\pux f)\ux+\Psi(f),\quad(\ux\,f)\pux =\ux(f \pux)+\Psi(f)$,
\item[(ii)] $\pux\Psi(f)=-2f\pux-\Psi(\pux f),\quad\Psi(f)\pux=-2\pux f-\Psi(f\pux)$,
\item[(iii)] $\pux\Psi(f)\pux=\Psi(\pux f\pux),\quad\pux^2\Psi(f)=\Psi(\pux^2 f)$.
\end{itemize}
\end{proposition}

\section{Auxiliary results}

Before stating the main theorems to be proved, we need some auxiliary results to which this section is devoted. We frequently use the symbol $\B(\ux,r)$ ($\partial\B(\ux,r)$) to denote the open ball (sphere) in $\R^m$ with center $\ux$ and radius $r$. In particular, $\B=\B(0,1)$. The volume (area) of $\B(\ux,r)$ ($\partial\B(\ux,r)$) will be denoted by $|\B(\ux,r)|$ ($|\partial\B(\ux,r)|$).
\begin{lemma}\label{fundamental}
Let $f$ be an $\R_{0,m}$-valued function, continuous in $\overline{\B}$, two-sided three-monogenic in $\B$ and vanishing on $\partial\B$. Then, there exists a two-sided monogenic function $\eta$ such that
\begin{equation}\label{representation}
f(\ux)=(1-\|\ux\|^2)\eta(\ux),\,\,\ux\in\B.
\end{equation}
\end{lemma}
\pf By the Almansi-type decomposition for three-monogenic functions established in \cite{MalRen}, there exists a two-sided monogenic function $\varphi$ and a harmonic function $h$ in $\B$ such that
\[
f(\ux)=\|\ux\|^2\varphi(\ux)+h(\ux),\,\,\ux\in\B.
\]
Since $f$ is uniformly continuous in $\overline\B$ and vanished on $\partial\B$, for any $\varepsilon>0$, there exists $r_{\varepsilon}>0$ such that $\|f(\ux)\|<\varepsilon$ for every $\ux$ in the annular region $\{r_{\varepsilon}<\|\ux\|<1\}$.

Let us introduce the sequence $\{r_n\}_{n\geq 1}$, with $r_{\varepsilon}=r_1<r_2<\cdots r_n<\cdots<1$, converging to $1$. 

Now, let $\ux_0$ be a fixed, but arbitrary, point in $\B$ and choose $N>0$ such that $\ux_0\in\B(0,r_n)$ for $n\geq N$.

Because $\|r_n^2\varphi(\xi_n)+h(\xi_n)\|=\|f(\xi_n)\|<\varepsilon$ for $\xi_n\in\partial\B(0,r_n)$, the harmonicity of the function $r_n^2\varphi+h$ together with the maximum modulus principle for harmonic functions yield
\[
\|r_n^2\varphi(\ux_0)+h(\ux_0)\|<\varepsilon.
\] 
Letting $n\to\infty$ we have $\|\varphi(\ux_0)+h(\ux_0)\|\le\varepsilon$, which implies $\varphi=-h$ and the proof is completed with $\eta=-\varphi$.\qed
\begin{lemma}\label{harmonic}
Let $\omega$ be an $\R$-valued harmonic function in a domain $\Omega$ containing the closed ball $\overline\B(\ux,r)\subset\R^m$. Then
\[
\left|\partial_{x_i}\omega(\ux)\right|\le \frac{m}{r}\sup_{\uy\in\B(\ux,r)}|\omega(\uy)-\omega(\ux)|.
\]
\end{lemma} 
\pf Since $\frac{\partial\omega}{\partial x_i}$ is harmonic in $\overline\B(\ux,r)$, the mean-value property leads to
\[
\partial_{x_i}\omega(\ux)=\frac{1}{|\B(\ux,r)|}\int\limits_{\B(\ux,r)}\partial_{y_i}\omega(\uy)d\uy,
\]
which, after using the Gauss formula, gives
\[
\partial_{x_i}\omega(\ux)=\frac{1}{|\B(\ux,r)|}\int\limits_{\partial\B(\ux,r)}\theta_i(\uy)\omega(\uy)d\uy,
\]
where $\theta_i$ denotes the $i$-th component of the unit normal vector at $\uy$. 

Consequently, we have
\[
\partial_{x_i}\omega(\ux)=\frac{1}{|\B(\ux,r)|}\int\limits_{\partial\B(\ux,r)}\theta_i(\uy)(\omega(\uy)-\omega(\ux))d\uy
\]
and so
\[
\left|\partial_{x_i}\omega(\ux)\right|\le\frac{|\partial\B(\ux,r)|}{|\B(\ux,r)|}\sup_{\uy\in\B(\ux,r)}|\omega(\uy)-\omega(\ux)|.
\]
Since $|\B(\ux,r)|=\frac{r}{m}|\partial\B(\ux,r)|$, the desired estimate follows.\qed

Using standard componentwise arguments, the above result easily extends to $\R_{0,m}$-valued harmonic functions if the Clifford norm $\|\cdot\|$ is used instead of modulus. 
\section{A solvability criterion in odd-dimensional spaces}
In this section we establish a necessary and sufficient condition for the solvability of the Dirichlet problem \eqref{DPI} in odd dimensions. Given a continuous function $g$ in $\partial\B$ we will consider its Poisson transform by
\[
\cP[g](\ux):=\frac{1}{|\partial\B|}\int\limits_{\partial\B}\frac{1-\|\ux\|^2}{\|\uy-\ux\|^m}g(\uy)d\uy.
\]
First, we need the following
\begin{proposition}\label{Prop infra representation}
Let $f$ be an $\R_{0,m}$-valued inframonogenic function in $\B\subset\R^m$. Suppose $f$ to be continuous in the closed ball $\overline\B$ with trace $g$ on $\partial\B$. Then, for each $\ux\in\B$ we have:
\begin{equation}\label{infra representation}
f(\ux)=\frac{1}{2}(1-\|\ux\|^2)\Psi^{-1}(\pux\cP[g]\pux)+\cP[g](\ux).
\end{equation}  
\end{proposition}
\pf Since $f$ is inframonogenic in $\B$, it is two-sided three-monogenic there. The same happens with the function $f-\cP[g]$, which is continuous in $\overline\B$ and has trace identically zero on $\partial\B$. By Lemma \ref{fundamental}, there exists a two-sided monogenic function $\eta$ such that 
\[
f(\ux)-\cP[g](\ux)=(1-\|\ux\|^2)\eta(\ux),\,\,\,\ux\in\B.
\] 
Applying the sandwich operator $\pux (\cdot)\pux$ in both hand-sides of the previous identity leads to
\[
\pux((1-\|\ux\|^2)\eta(\ux))\pux=-\pux(\cP[g](\ux))\pux.
\]
From the direct identity $\pux((1-\|\ux\|^2)\eta)=-2\ux\eta$ and  Proposition \ref{p1}-$(1)$, it follows that 
\[\pux((1-\|\ux\|^2)\eta)\pux=-2\Psi(\eta),\]
where use has been made of the two-sided monogenicity of $\eta$.

Consequently,  we have $$\Psi(\eta)=\frac{1}{2}\pux(\cP[g](\ux))\pux,$$ which together with the invertibility of $\Psi$ yield  $\eta=\frac{1}{2}\Psi^{-1}(\pux(\cP[g])\pux)$ and we are done.\qed
\begin{theorem}[Solvability Criterion]\label{criterion}
Let $g\in C(\partial\B)$, then the Dirichlet problem \eqref{DPI} is solvable if and only if
\begin{equation}\label{condition}
\lim_{\|\ux\|\to 1-} (1-\|\ux\|)\pux(\cP[g])\pux=0.
\end{equation}
\end{theorem}
\pf If $f$ is a solution of \eqref{DPI}, then by Proposition \ref{Prop infra representation} we have
\[
f(\ux)=\frac{1}{2}(1-\|\ux\|^2)\Psi^{-1}(\pux\cP[g]\pux)+\cP[g](\ux).
\]
Then
\[
\lim_{\|\ux\|\to 1-} (1-\|\ux\|)\Psi^{-1}(\pux(\cP[g])\pux)=\frac{1}{2}\lim_{\|\ux\|\to 1-} (1-\|\ux\|^2)\Psi^{-1}(\pux(\cP[g])\pux)=\lim_{\|\ux\|\to 1-} (f(\ux)-\cP[g](\ux))=0.
\]
By linearity we have 
\[
\lim_{\|\ux\|\to 1-} \Psi^{-1}((1-\|\ux\|)\pux(\cP[g])\pux)=\lim_{\|\ux\|\to 1-} (1-\|\ux\|)\Psi^{-1}(\pux(\cP[g])\pux)=0
\]
and hence
\[
\lim_{\|\ux\|\to 1-} (1-\|\ux\|)\pux(\cP[g])\pux=0.
\]
Conversely, if \eqref{condition} is fulfilled then the function
\[
f(\ux)=\frac{1}{2}(1-\|\ux\|^2)\Psi^{-1}(\pux\cP[g]\pux)+\cP[g](\ux).
\]
is obviously a solution of \eqref{DPI}.\qed
The following result about the uniqueness of the solutions of the problem \eqref{DPI} was already established in \cite{M5}. We include here a simpler proof, which made essential use of Lemma \ref{fundamental}. 
\begin{theorem}[Uniqueness]\label{odduniqueness}
If $f\in C^2(\B)\cap C(\overline\B)$ is inframonogenic in $\B$ and has null trace on $\partial\B$, then $f\equiv 0$ in $\B$.
\end{theorem}
\pf By Lemma \ref{fundamental}, there exists a two-sided monogenic function $\eta$ in $\B$ such that $f=(1-\|\ux\|^2)\eta(\ux)$. 

Applying the sandwich operator and repeating the calculations carried out in the proof of Proposition \ref{Prop infra representation}, we obtain that $\Psi(\eta)=0$ in $\B$ and  so $\eta=0$ there. Hence $f=(1-\|\ux\|^2)\eta(\ux)=0$ in $\B$.\qed 
\subsection{Boundary data in $C^{1}(\partial\B)$}
Here we use the above results in the case of data with some degree of smoothness on the boundary.  Let us begin with the following 
\begin{theorem}\label{caso c1}
Let $h\in C^1(\overline\B)$ be harmonic in $\B$. Then
\[
\lim_{\|\ux\|\to 1-} (1-\|\ux\|)\partial_{x_i}\partial_{x_j}h=0.
\] 
\end{theorem}
\pf Let $\omega_j=\partial_{x_j}h$ and choose an arbitrary point $\ux$ in $\B$. For $\delta=\frac{1-\|\ux\|}{2}$, we have the inclusion $\overline{\B(\ux,\delta)}\subset\B$.

Since $\omega_j$ is harmonic in $\B$, by Lemma \ref{harmonic} we have:
\[
\|\partial_{x_i}\omega_j(\ux)\|\le\frac{m}{\delta}\sup_{\uy\in\B(\ux,\delta)}\|\omega_j(\uy)-\omega_j(\ux)\|,
\]
which after multiplying by $1-\|\ux\|$ leads to
\[
\|(1-\|\ux\|)\partial_{x_i}\omega_j(\ux)\|\le 2\delta\frac{m}{\delta}\sup_{\uy\in\B(\ux,\delta)}\|\omega_j(\uy)-\omega_j(\ux)\|=2m\sup_{\uy\in\B(\ux,\delta)}\|\omega_j(\uy)-\omega_j(\ux)\|.
\]
If $\|\ux\|\to 1-$, then $\delta\to 0$ and the right-hand side of the above inequality goes to $0$ in virtue of the uniform continuity of $\omega_j$ in $\overline\B$. The proof has been completed.\qed
\begin{corollary}\label{coroc1}
Let $g\in C^1(\partial\B)$. Then, the Dirichlet problem \eqref{DPI} has a unique solution, which is given by \eqref{infra representation}. 
\end{corollary}
\pf The uniqueness has been already proved. Moreover, by the classical regularity theorem for the Poisson integral, it follows that $\cP[g]$ belongs to $C^1(\overline\B)$ as well. Then, the existence follows readily by expanding $(1-\|\ux\|)\pux(\cP[g])\pux$ into its real components, as a direct consequence of Theorem \ref{criterion} and Theorem \ref{caso c1}.\qed

It should be noted that the above theorem gives a positive answer to the open question left in \cite[p.507]{M5}.

We now turn to the question of whether $C^1(\partial\B)$ is the optimal class in which the Dirichlet problem \eqref{DPI} admits a solution. Indeed, we shall prove that this condition is sharp, in the sense that one can construct a continuous function which is not of class $C^1(\partial\B)$, for which the problem has no solution. 

To this end, we begin with the following
\begin{lemma}\label{contraejemplo}
There exist an $\R$-valued harmonic function in $\B$, continuous in $\overline\B$ and such that
\[
\lim_{\|\ux\|\to 1-} (1-\|\ux\|)\partial_{x_i}\partial_{x_j}h\not=0,
\] 
for some $i,j\in\{1,2,\dots, m\}$.
\end{lemma}
\pf Consider the complex valued function
\[
F(z)=\sum\limits_{n=1}^\infty\frac{z^n}{n^2},\,\,|z|\le 1.
\]
Obviously, the above series converges absolutely and uniformly in the closed unit disk $\overline\cD$ and so it is continuous in $\overline\cD$ and holomorphic in $\cD$. Moreover, we have in $\cD$:
\[
F'(z)=\sum\limits_{n=1}^\infty\frac{z^{n-1}}{n},\,\,\,F''(z)=\sum\limits_{n=2}^\infty\frac{n-1}{n}z^{n-2}.
\]
It should be noticed that $F'(z)$ diverges for $z=1$ and so $F\notin C^1(\overline\cD)$.

By using the Taylor expansion of the principal branch of the logarithm (Newton-Mercator series), we have
\[
-\log(1-z)=\sum\limits_{n=1}^\infty\frac{z^n}{n},\,\,|z|< 1,
\]  
which yields
\[
F'(z)=\frac{1}{z}\sum\limits_{n=1}^\infty\frac{z^n}{n}=-\frac{1}{z}\log(1-z),\,\,0<|z|< 1,
\]
and
\[
F''(z)=\frac{1}{z(1-z)}+\frac{\log(1-z)}{z^2},\,\,0<|z|< 1.
\]
If we restrict $z$ to the open interval $(-1,1)$, we have
\[(1-t)F''(t)=\frac{1}{t}+\frac{(1-t)\log(1-t)}{t^2},\,\,\, t\in (0,1),\]
and consequently $\lim\limits_{t\to 1-}(1-t)F''(t)=1$.

Let us now define in $\B$ the $\R$-valued harmonic function
\[
h(\ux)=h(x_1,x_2,\dots,x_m)=\Re(F(x_1+ix_2)),
\]
where as usual $\Re (z)$ denotes the real part of the complex number $z$.

Direct computation gives
\[
\partial_{x_1}^2 h(t,0,0\dots,0)=F''(t).
\]
Thus, taking the path determined by the real axis of the first component, we have
\begin{equation}\label{contra}
\lim\limits_{x_1\to 1-}(1-\|\ux\|)\partial_{x_1}^2h(\ux)=\lim\limits_{t\to 1-}(1-t)F''(t)=1,
\end{equation}
which completes the proof.\qed

Now let $g(\ux)=h(\ux)e_1$, where $h$ is the function constructed in Lemma \ref{contraejemplo}. Of course, we have that $g\in C(\partial\B)$. Then, by Theorem \ref{criterion}, the problem \eqref{DPI} has solution if and only if  
\[
\lim_{\|\ux\|\to 1-} (1-\|\ux\|)\pux(\cP[g])\pux=0.
\]
Since $g$ is harmonic by construction, we have $g=\cP[g]$ in $\B$. On the other hand,
\[
\pux g(\ux)\pux=\sum\limits_{i,j\in\{1,2,\dots,m\}}e_i\partial_{x_i}\partial_{x_j}g(\ux)e_j=e_1\partial_{x_1}^2h(\ux)e_1e_1+e_1\partial_{x_1}\partial_{x_2}h(\ux)e_1e_2+e_2\partial_{x_2}^2h(\ux)e_1e_2+e_2\partial_{x_2}\partial_{x_1}h(\ux)e_1e_1.
\]
After some simplification and using the fact that  $h$ is harmonic, we obtain
\[
\pux g(\ux)\pux=-2\partial_{x_1}^2h(\ux)e_1-2\partial_{x_1}\partial_{x_2}h(\ux)e_2
\]
and hence
\[
(1-\|\ux\|)\pux(\cP[g])\pux=-2(1-\|\ux\|)\partial_{x_1}^2h(\ux)e_1-2(1-\|\ux\|)\partial_{x_1}\partial_{x_2}h(\ux)e_2.
\]
The last expression goes to $0$ as $\|\ux\|\to 1-$, if and only if each one of the summands in the right-hand side does so, which does no hold in virtue of \eqref{contra}.
\subsection{Boundary data in $C^{1,\nu}(\partial\B)$}
It is a well-known fact that the solution of the Dirichlet problem for harmonic functions preserves the regularity of the boundary data. In this subsection, we present a weak form of the regularity theorem for the Dirichlet problem \eqref{DPI}. Indeed, we prove that the solution of \eqref{DPI} is in $C^\nu(\overline\B)$, whenever $g\in C^{1,\nu}(\partial\B)$. In other words, the solution of the Dirichlet problem reduces the smoothness of the boundary data by one degree.
\begin{lemma}\label{L1nu}
Let $h$ be an $\R$-valued harmonic function in $\B$ and $h\in C^{1,\nu}(\overline\B)$ ($0<\nu<1$). Then, there exists a constant $C$ such that
\[
|\partial_{x_i}\partial_{x_j}h(\ux)|<C(1-\|\ux\|)^{\nu-1},\,\,\ux\in\B.
\]
\end{lemma}
\pf Let $\ux\in\B$ and consider the ball $\B(\ux,r)$, with $r=\frac{1}{2}(1-\|\ux\|)$. By Lemma \ref{harmonic} applied to the function $\partial_{x_j} h$, we have
\[
|\partial_{ x_i}\partial_{x_j}h(\ux)|\le\frac{m}{r}\sup\limits_{\uy\in\B(\ux,r)}|\partial_{x_j}h(\uy)-\partial_{x_j}h(\ux)|\le C\frac{m}{r}\|\ux-\uy\|^\nu\le C m r^{\nu-1}\le C (1-\|\ux\|)^{\nu-1}.
\]
\qed
Throughout, $C$ denotes a generic constant, possibly changing from line to line.
\begin{lemma}\label{L2nu}
Under the same assumptions of Lemma \ref{L1nu}, there exists a constant $C$ such that
\[
|\partial_{x_k}\partial_{x_i}\partial_{x_j}h(\ux)|<C(1-\|\ux\|)^{\nu-2},\,\,\ux\in\B.
\] 
\end{lemma}
\pf Let $\ux\in\B$ and once again consider the ball $\B(\ux,r)$, with radius $r=\frac{1}{2}(1-\|\ux\|)$. Lemma \ref{harmonic} applied to $\partial_{x_i}\partial_{x_j}h$, gives
\begin{equation}\label{eq1nu}
|\partial_{x_k}\partial_{ x_i}\partial_{x_j}h(\ux)|\le\frac{m}{r}\sup\limits_{\uy\in\B(\ux,r)}|\partial_{x_i}\partial_{x_j}h(\uy)-\partial_{x_i}\partial_{x_j}h(\ux)|\le \frac{C}{r}((1-\|\uy\|)^{\nu-1}+(1-\|\ux\|)^{\nu-1}),
\end{equation}
where use has been made of Lemma \ref{L1nu}.

Moreover, one can readily verify that $1-\|\ux\|\le 2(1-\|\uy\|)$ for $\uy\in\B(\ux,r)$. Consequently, we have 
\[(1-\|\uy\|)^{\nu-1}\le \frac{1}{2^{\nu-1}}(1-\|\ux\|)^{\nu-1},\] 
which together with \eqref{eq1nu} complete the proof.\qed  
We are now in a position to state and prove the following
\begin{theorem}\label{auxiliar1nu}
Let $h$ be an $\R$-valued harmonic function in $\B$ and $h\in C^{1,\nu}(\overline\B)$ ($0<\nu<1$). Then, the function 
\[F(\ux)=(1-\|\ux\|^2)\partial_{x_i}\partial_{x_j}h\]
belongs to $C^{0,\nu}(\overline\B)$.
\end{theorem}
\pf We need to prove that there exists a positive constant $C$ such that:
\[
|F(\ux)-F(\uy)|\le C\|\ux-\uy\|^\nu,\,\,\ux,\uy\in\overline\B.
\]  
We will consider two cases. Namely, the case when $\|\ux-\uy\|\geq\frac{1}{2}(1-\|\ux\|)$ and the case when $\|\ux-\uy\|<\frac{1}{2}(1-\|\ux\|)$. Without lost of generality we can assume that $\|\ux\|\geq\|\uy\|$.

In the first case, we have $1-\|\uy\|\le 1-\|\ux\|+\|\ux-\uy\|\le 3\|\ux-\uy\|$.

Therefore
\[
|F(\ux)-F(\uy)|\le |F(\ux)|+|F(\uy)|=(1-\|\ux\|^2)|\partial_{x_i}\partial_{x_j}h(\ux)|+(1-\|\uy\|^2)|\partial_{x_i}\partial_{x_j}h(\uy)|,
\]
which after using Lemma \ref{L1nu} yields
\[
|F(\ux)-F(\uy)|\le C[(1-\|\ux\|)^\nu+(1-\|\uy\|)^\nu]\le C(2^\nu+3^\nu)\|\ux-\uy\|^\nu.
\]
Let us now consider the second case, namely when $\|\ux-\uy\|<\frac{1}{2}(1-\|\ux\|)$.

First notice that
\[
\partial_{x_k}F(\ux)=-2x_k\partial_{x_i}\partial_{x_j}h(\ux)+(1-\|\ux\|^2)|\partial_{x_k}\partial_{x_i}\partial_{x_j}h(\ux)
\]
and consequently, it follows from Lemma \ref{L1nu} and Lemma \ref{L2nu} that
\[
|\partial_{x_k}F(\ux)|\le 2C(1-\|\ux\|)^{\nu-1}+2C (1-\|\ux\|)(1-\|\ux\|)^{\nu-2}\le C(1-\|\ux\|)^{\nu-1}
\]
and hence
\begin{equation}\label{gradient}
\|\pux F(\ux)\|\le C(1-\|\ux\|)^{\nu-1}.
\end{equation}
The Mean Value Theorem applied to $F$ in the convex set $\overline\B$ ensures that
\[
|F(\ux)-F(\uy)|\le \|\pux F(\xi)\|\|\ux-\uy\|,
\]
for certain point $\xi$ in the segment joining $\ux$ and $\uy$.

The latter together with \eqref{gradient} yields
\[
|F(\ux)-F(\uy)|\le C(1-\|\xi\|)^{\nu-1}\|\ux-\uy\|.
\]
Finally, taking into account that $\xi$ lies in the segment joining $\ux$ and $\uy$ and the fact that $\|\ux-\uy\|<\frac{1}{2}(1-\|\ux\|)$ we have
$1-\|\xi\|\geq \frac{1}{2}(1-\|\ux\|)$ and so
\[
(1-\|\xi\|)^{\nu-1}\le\frac{1}{2^{\nu-1}}(1-\|\ux\|)^{\nu-1}\le \frac{2^{\nu-1}}{2^{\nu-1}}\bigg(\frac{1-\|\ux\|}{2}\bigg)^{\nu-1}\le \|\ux-\uy\|^{\nu-1}.
\]
The proof is completed.\qed

As a corollary, we now state the main result of this section.
\begin{corollary}\label{final1nu}
Let $g\in C^{1,\nu}(\partial\B)$. Then, the solution of the Dirichlet problem \eqref{DPI} belongs to $C^{0,\nu}(\overline\B)$. 
\end{corollary}
\pf The standard regularity result says that the Poisson transform $\cP[g]$ preserves the smoothness of $g$. Namely, $\cP[g]\in C^{1,\nu}(\overline\B)$. With this at hand, the statement follows directly from Corollary \ref{coroc1} and Theorem \ref{auxiliar1nu}.\qed
\subsection{Boundary data in $C^{k,\nu}(\partial\B)$}
Here we extend the previous section to arbitrary degree of smoothness. 
\begin{theorem}
Let $h$ be an $\R$-valued harmonic function in $\B$ and $h\in C^{k,\nu}(\overline\B)$ ($k\geq 1$). Then, the function 
\[F(\ux)=(1-\|\ux\|^2)\partial_{x_i}\partial_{x_j}h\]
belongs to $C^{k-1,\nu}(\overline\B)$.
\end{theorem}

\pf
The proof is carried out by induction on $k$. The case $k=1$ was already established in Theorem \ref{auxiliar1nu}. Assume the result holds for some $k\geq 1$ and prove it for $k+1$.  

Let $h\in C^{k+1,\nu}(\overline\B)$ and let $\beta=(\beta_1,\beta_2,\dots,\beta_n)$ a multiindex of degree $|\beta|=\beta_1+\beta_2+\cdots +\beta_n=k$. It is readily seen that
\[
\partial^\beta=\partial_{x_l}\partial^\gamma
\]  
for some variable $x_l$ and some multiindex $\gamma$ with $|\gamma|=k-1$.

Therefore, we have
\begin{equation}\label{induction}
\partial^{\beta}F=\partial^{\gamma}\partial_{x_l}\bigg[(1-\|\ux\|^2)\partial_{x_i}\partial_{x_j}h\bigg]=\partial^\gamma\bigg[-2x_l\partial_{x_i}\partial_{x_j}h+(1-\|\ux\|^2)\partial_{x_i}\partial_{x_j}(\partial_{x_l}h)\bigg]=\partial^\gamma\bigg[-2x_l\partial_{x_i}\partial_{x_j}h\bigg]+\partial^\gamma\bigg[(1-\|\ux\|^2)\partial_{x_i}\partial_{x_j}(\partial_{x_l}h)\bigg].
\end{equation} 
Since $h\in C^{k+1,\nu}(\overline\B)$ and $|\gamma|=k-1$, the first summand in \eqref{induction} belongs to $C^{0,\nu}(\overline\B)$. On the other hand, since $\partial_{x_l}h\in C^{k,\nu}(\overline\B)$, the induction hypothesis ensures that the function $(1-\|\ux\|^2)\partial_{x_i}\partial_{x_j}(\partial_{x_l}h)$ belongs to $C^{k-1,\nu}(\overline\B)$ and hence the second summand in \eqref{induction} belongs to $C^{0,\nu}(\overline\B)$, as well. The proof is completed. \qed

The following is a direct consequence of the previous result.
\begin{corollary}\label{finalknu}
Let $g\in C^{k,\nu}(\partial\B)$. Then, the solution of the Dirichlet problem \eqref{DPI} belongs to $C^{k-1,\nu}(\overline\B)$. 
\end{corollary}

We consider an alternative perspective on our results. For a given $\R_{0,m}$-valued function defined in the sphere $\partial\B$ of $\R^m$ (with $m$ odd) let us introduce the inframonogenic Poisson transform
\[
\cP_{\I}[g]=\frac{1}{2}(1-\|\ux\|^2)\Psi^{-1}(\pux\cP[g]\pux)+\cP[g](\ux),\,\,\ux\in\B.
\]
Then we can rephrase our achievements in the following
\begin{theorem}
For $k\geq 1$, the inclusion
\[
\cP_{\I}\bigg(C^{k,\nu}(\partial\B)\bigg)\subset C^{k-1,\nu}(\overline\B)
\]
holds.
\end{theorem}
\section{Even dimensional case}
The present section will be brief. From the beginning, the significance of the operator $\Psi$  was evident: in odd dimensions it works as a linear isomorphism, whereas in even dimensions this fails. Nevertheless, the method of the previous subsections may be applied with success to the solution of \eqref{DPI} in any dimension $m$,
except when $m =2p$ and the $p$-vector part $[g]_p$ of the boundary data $g$ is not zero. In this case the operator $\Psi^{-1}$ is no longer available and the method breaks down. In \cite{M5} it was proved how, in this situation the Dirichlet problem \eqref{DPI} may be unsolvable even when the boundary condition is
given by a polynomial. The proof of this result relied on a necessary an sufficient condition for the solvability established in \cite[Theorem 6]{M5}. Here we present an improved version (the solution does not need to be smooth in $\overline\B$) with a simpler proof based on Lemma \ref{fundamental}.
\begin{theorem}\label{T6}
Let $p=\frac{m}{2}$ and let $g_p$ be a $p$-vector valued function in $C(\partial{\B})$. Then, the Dirichlet problem \eqref{DPI} with $g=g_p$ has solution in $C^2(\B)\cap C(\overline{\B})$ if and only if
\begin{equation}\label{Condition}
\pux(\cP[g_p])\pux=0.
\end{equation} 
\end{theorem}
\pf If $f$ is a solution of \eqref{DPI}, then by Lemma \ref{fundamental} there exists a two-sided monogenic function $\eta$, this time $p$-vector valued, such that 
\[
f(\ux)-\cP[g_p](\ux)=(1-\|\ux\|^2)\eta(\ux),\,\,\,\ux\in\B.
\] 
Applying the sandwich operator $\pux (\cdot)\pux$ to both sides of the previous identity yields
\[
\pux((1-\|\ux\|^2)\eta)\pux=-\pux(\cP[g_p])\pux.
\]
As proved earlier, we have $\pux((1-\|\ux\|^2)\eta)\pux=-2\Psi(\eta)$.  

Consequently,  we have $\Psi(\eta)=\frac{1}{2}\pux(\cP[g_p])\pux$. Since $\eta$ is $p$-vector valued, it follows from \eqref{=0} that $\Psi(\eta)=0$ and hence $\pux(\cP[g_p])\pux=0$. 

The converse is immediate.\qed 

The above procedure shows that in even dimension there exist infinitely many functions of the form $(1-\|\ux\|^2)\eta_p$, where  $\eta_p$ is $p$-vector valued monogenic in $\B$, with $p=\frac{m}{2}$, which solve the homogeneous Dirichlet problem
\begin{equation}\label{DPI0}
\Biggl\{
\begin{array}{rl}
\pux f\pux=0 &\,\,\mbox{in}\,\,\B,\\f=0&\,\,\mbox{in}\,\,\partial\B.
\end{array} 
\end{equation}
Consequently, no uniqueness theorem holds in this even-dimensional case. 
%%%%%%%%%%%%%%%%%%%%%%%%%%%%%%%%%%%%%%%%%%%%%%%%%%%%%%%
\section{A generalized framework}  
In order to derive explicit representations of the isomorphisms constructed in \cite{das1}, we shall restrict our analysis in this section to the three-dimensional setting $\R^3\cong\R_{0,3}^{(1)}$.

Let $\chi=\{\chi_1,\chi_2,\chi_3\}\subset\R_{0,3}^{(1)}$. 
By a slight abuse of notation, we shall write
$-\chi:=\{-\chi_1,-\chi_2,-\chi_3\}$.
For a domain $\Omega\subset\R^3$, we define on $C^1(\Omega)$ the following
non-standard Dirac operator, which will be referred to as the
$\chi$-Dirac operator:
\begin{equation*}
\hD:=\chi_1\partial_{x_1}+
\chi_2\partial_{x_2}+
\chi_3\partial_{x_3}.
\end{equation*}

A direct computation shows that the factorization
\begin{equation}\label{Fac Lap}
\hD\hD=-\Delta
\end{equation}
holds if and only if
\begin{equation}\label{Estr Cond}
\chi_i\chi_j+\chi_j\chi_i
=-2\delta_{ij},
\qquad i,j=1,2,3,
\end{equation}
where $\delta_{ij}$ denotes the Kronecker delta. Consequently, the
factorization \eqref{Fac Lap} holds precisely when $\chi$ constitutes an
orthonormal basis of
$\R^3$.

Any set satisfying \eqref{Estr Cond} is called a
\emph{structural set}. The canonical example is the standard structural set
\[
\chi_{\rm st}:=\{e_1,e_2,e_3\},
\]
which gives rise to the standard Dirac operator
\[
\pux\equiv{}^{\chi_{\rm st}}\!\pux.
\]

The use of structural sets provides a flexible framework for extending and
revisiting several topics in Clifford analysis from a broader
geometric perspective. In particular, structural sets have been considered
in connection with the mapping properties of multidimensional
Ahlfors--Beurling transforms (or $\Pi$-operators), $\overline{\partial}$-problems,
generalized Kolosov--Muskhelishvili formulas, $M$-conformal mappings, and
additive decompositions of contragenic polynomials. They also provide useful
tools for studying the composition of monogenic functions with M\"obius
transformations and the reciprocals of monogenic functions; see
\cite{ss1,ss2,ss3,nguyen,ss4}.

We now consider a class of second-order partial differential equations
associated with two structural sets
\[
\chi=\{\chi_1,\chi_2,\chi_3\},
\qquad
\vartheta=\{\vartheta_1,\vartheta_2,\vartheta_3\}.
\]
More precisely, we study the general sandwich equation
\begin{equation}\label{sge}
\hD f\pD=0,
\end{equation}
whose solutions are referred to as
$(\chi,\vartheta)$-inframonogenic functions; see \cite{das}.

The remainder of this section is devoted to the Dirichlet problem
\begin{equation}\label{DPIG}
\left\{
\begin{array}{rll}
\hD f\pD&=0, & \text{in }\B,\\[1mm]
f&=g, & \text{in }\partial\B.
\end{array}
\right.
\end{equation}

The algebraic isomorphisms developed in \cite{das1} play a central role in
our approach. They allow us to reduce the general sandwich equation
\eqref{sge} to the standard sandwich equation and, in particular, to
obtain explicit representations for the solution of the corresponding
Dirichlet problem. More precisely, one can construct an invertible Clifford
number $\varpi$ satisfying
\[
\varpi\chi_i=\vartheta_i\varpi,
\qquad i=1,2,3.
\]
Consequently, we obtain the equivalence
\[
\varpi\hD f\pD=0
\quad\Longleftrightarrow\quad
\pD\varpi f\pD=0.
\]

Let
\[
\mathcal{M}_{\vartheta,\chi}
=
\begin{pmatrix}
c_{11}&c_{21}&c_{31}\\
c_{12}&c_{22}&c_{32}\\
c_{13}&c_{23}&c_{33}
\end{pmatrix}
\]
be the change-of-basis matrix from $\vartheta$ to $\chi$. If
$\mathcal{M}_{\vartheta,\chi}$ is non-symmetric and has unit determinant, then the
invertible Clifford number
\[
\varpi=
-1-c_{11}-c_{22}-c_{33}
+(c_{12}-c_{21})\vartheta_1\vartheta_2
+(c_{13}-c_{31})\vartheta_1\vartheta_3
+(c_{23}-c_{32})\vartheta_2\vartheta_3
\]
provides the required intertwining relation. We next make use of the correspondence
\[
\cS_{\vartheta}:
\sum_Aa_A\vartheta_A
\longmapsto
\sum_Aa_Ae_A.
\]
Applying $\cS_{\vartheta}$ to the transformed equation gives
\[
\pux\cS_{\vartheta}[\varpi f]\pux=0,
\]
which is the standard sandwich equation. Hence, the original problem
can be reduced to a Dirichlet problem for usual inframonogenic functions.

When $\mathcal{M}_{\vartheta,\chi}$ is symmetric, three principal cases have to be
considered, namely $c_{13}\neq0$, $c_{23}\neq0$, or $c_{12}\neq0$. In the
first two cases, we may choose
\[
\varpi
=
(-c_{33}-1)\vartheta_1\vartheta_2
+c_{23}\vartheta_1\vartheta_3
-c_{13}\vartheta_2\vartheta_3,
\]
whereas in the remaining case we take
\[
\varpi
=
-c_{23}\vartheta_1\vartheta_2
+(c_{22}+1)\vartheta_1\vartheta_3
-c_{12}\vartheta_2\vartheta_3.
\]
The remaining configurations are exceptional cases that can be treated
separately. In each of them, an explicit inver\-tible Clifford number
$\varpi$ can still be constructed.

Finally, if $\mathcal{M}_{\vartheta,\chi}$ has determinant $-1$, we use the
equivalence between
\[
\hD f\pD=0
\qquad\text{and}\qquad
{}^{-\chi}\!\pux\,f\pD=0,
\]
which have the same solution space. Thus, replacing $\chi$ by
$-\chi$, we may reduce the problem to the case in which
$\mathcal{M}_{\vartheta,-\chi}$ has unit determinant.

The preceding analysis, together with Corollaries \ref{coroc1} and \ref{final1nu}, yields the
following result.

\begin{theorem}\label{TDPG}
Let $\chi$ and $\vartheta$ be arbitrary structural sets in $\mathbb{R}^3$, and let
$g\in C^{1,\nu}(\partial\mathbb{B})$ be given. Then the Dirichlet problem
\eqref{DPIG} admits a unique solution in $C^{0,\nu}(\overline{\mathbb{B}})$,
which is given explicitly by the formula
\begin{equation}
f(\ux)=\varpi^{-1}\mathcal{S}_\vartheta^{-1}\bigg[\frac{1}{2}(1-\|\ux\|^2)\Psi^{-1}(\pux\cP[\mathcal{S}_\vartheta\varpi g]\pux)+\cP[\mathcal{S}_\vartheta\varpi g](\ux)\bigg]. 
\end{equation}

\end{theorem}

We conclude with an explicit example illustrating the reduction procedure
described above. A similar example can be found in \cite{M5}; however, we
include it here for completeness and to make the construction transparent.

\begin{example}
Consider the structural sets
\begin{align*}
\chi
&=
\left\{
\frac{\sqrt{2}}{2}e_1+\frac{\sqrt{2}}{2}e_2,\,
\frac{\sqrt{2}}{2}e_1-\frac{\sqrt{2}}{2}e_2,\,
-e_3
\right\},\\
\vartheta
&=
\left\{
\frac{\sqrt{3}}{3}e_1-\frac{\sqrt{2}}{2}e_2-\frac{\sqrt{6}}{6}e_3,\,
\frac{\sqrt{3}}{3}e_1+\frac{\sqrt{2}}{2}e_2-\frac{\sqrt{6}}{6}e_3,\,
\frac{\sqrt{3}}{3}e_1+\frac{\sqrt{6}}{3}e_3
\right\}.
\end{align*}

We consider the Dirichlet problem
\begin{equation}\label{PGE}
\left\{
\begin{array}{rll}
\hD f\pD&=0, & \text{in }\B,\\[1mm]
f&=x_3^2e_2, & \text{in }\partial\B.
\end{array}
\right.
\end{equation}
For these structural sets, the Clifford number
\[
\varpi=
\left(\frac{\sqrt{6}}{3}-1\right)
\vartheta_1\vartheta_2
+\frac{\sqrt{6}}{6}\vartheta_1\vartheta_3
-\frac{\sqrt{6}}{6}\vartheta_2\vartheta_3
=
\left(1-\frac{\sqrt{6}}{3}\right)e_1e_2
-\frac{\sqrt{3}}{3}e_2e_3
\]
is invertible and provides the required transformation. Hence,
\eqref{PGE} is equivalent to
\begin{equation}\label{PGEt1}
\left\{
\begin{array}{rll}
\pD\varpi f\pD&=0, & \text{in }\B,\\[1mm]
\varpi f
&=
\left(\frac{\sqrt{6}}{3}-1\right)x_3^2e_1
-\frac{\sqrt{3}}{3}x_3^2e_3,
& \text{in }\partial\B.
\end{array}
\right.
\end{equation}
Applying $\cS_{\vartheta}$ to \eqref{PGEt1}, we obtain
\begin{equation}\label{Dps}
\left\{
\begin{array}{rll}
\pux\cS_{\vartheta}\varpi f\pux
&=0, & \text{in }\B,\\[1mm]
\cS_{\vartheta}\varpi f
&=
\left(\frac{\sqrt{2}}{2}-\frac{\sqrt{3}}{3}\right)x_3^2e_1
+\left(\frac{\sqrt{2}}{2}-\frac{\sqrt{3}}{3}\right)x_3^2e_2
-\frac{\sqrt{3}}{3}x_3^2e_3,
& \text{in }\partial\B.
\end{array}
\right.
\end{equation}
Thus, \eqref{Dps} is a Dirichlet problem for standard inframonogenic
functions. Its unique solution is
\[
\cS_{\vartheta}\varpi f
=
\left(\frac{\sqrt{3}}{3}-\frac{\sqrt{2}}{2}\right)
(x_1^2+x_2^2-1)(e_1+e_2)
-\frac{\sqrt{3}}{3}
(x_1^2+x_2^2+2x_3^2-1)e_3.
\]
Therefore,
\begin{align*}
f(\ux)
&=
\varpi^{-1}\cS_{\vartheta}^{-1}
\Big[
\left(\frac{\sqrt{3}}{3}-\frac{\sqrt{2}}{2}\right)
(x_1^2+x_2^2-1)(e_1+e_2)-\frac{\sqrt{3}}{3}
(x_1^2+x_2^2+2x_3^2-1)e_3
\Big].
\end{align*}
After simplifying, we obtain the explicit solution
\begin{align}\label{finalsol}
f(\ux)
&=
\left[
\frac{\sqrt{6}}{3}(\|\ux\|^2-1)+x_3^2
\right]e_2
-\frac{\sqrt{3}}{3}
(\|\ux\|^2-1)e_1e_2e_3.
\end{align}
It is worth emphasizing that, although the prescribed boundary data in
\eqref{PGE} is vector-valued, the corresponding solution
\eqref{finalsol} is not, in general, vector-valued. This is due to the
fact that the generalized sandwich operator
\[
\hD(\cdot)\pD
\]
does not, in general, preserve the space of $k$-vector-valued functions.
This is in contrast to the standard sandwich operator
\[
\pux(\cdot)\pux,
\]
which preserves the corresponding grade structure. Thus, this example demonstrates the effectiveness of the algebraic reduction
developed above and illustrates the broader structure arising in the general
setting. Moreover, many other interesting
examples can be constructed explicitly using the same approach.
\end{example}

A remarkable feature of this non-standard approach is the flexibility
provided by the choice of structural sets. In contrast to the standard
setting, where the differential operators are constructed from a fixed
orthonormal basis of the underlying Euclidean space, the use of arbitrary
structural sets allows one to adapt the associated sandwich operators to
different algebraic and analytic structures. In particular, given structural
sets $\chi$ and $\vartheta$, the operator
\[
\hD(\cdot)\pD
\]
generates a broad class of systems of partial differential equations whose properties
depend on the interaction between the two structural sets. This freedom makes
it possible to formulate and study general systems beyond
the standard operator $\pux(\cdot)\pux$.

This flexibility is particularly remarkable in the study of boundary value
problems. Different choices of the structural sets $\chi$ and $\vartheta$
give rise to different realizations of the underlying sandwich equation and,
consequently, to distinct classes of associated boundary value problems. The
choice of structural sets determines how the differential operators interact
with the Clifford algebra and may lead to different couplings among the
various grades of the unknown function, as well as to different forms of the
corresponding fundamental solutions and boundary operators. Thus, the
structural sets provide a natural degree of freedom that allows the formulation
of sandwich-type equations and their boundary value problems to be adapted to
the specific algebraic features of the problem under consideration.
\section*{Declarations}
\subsection*{Authors contributions}
The authors contributed equally to the manuscript and typed, read, and approved the final form of the manuscript, which
is the result of an intensive collaboration.
\subsection*{Conflict of interest statement}
The authors declare that they have no conflict of interest regarding the publication of this paper.

\subsection*{Data availability}
No datasets were generated or analyzed during the current study.

\subsection*{ORCID}
\noindent Arsenio Moreno Garc\'ia: \url{https://orcid.org/0000-0001-5984-5081}\\
Daniel Alfonso Santiesteban: \url{https://orcid.org/0000-0003-0248-3942}\\
Ricardo Abreu Blaya: \url{https://orcid.org/0000-0003-1453-7223}
%\section*{References}

\end{document}